\documentclass[a4paper,12pt]{amsart}

\usepackage{amsmath,amssymb}
\usepackage{amsbsy,amsfonts,array}
\usepackage{verbatim,graphicx}
\usepackage{epstopdf}
\usepackage{subfig}
\graphicspath{{./img/},{./}}
\usepackage[margin=3cm]{geometry} 
\usepackage[P,R]{dblfont}
\usepackage{hyperref}
\usepackage{moreverb}
\usepackage{graphicx}
\usepackage{multirow}
\usepackage{color}
\usepackage{psfrag}
\usepackage{stmaryrd}
\usepackage{epsfig}
\usepackage{bbold}
\usepackage{url}

\newcommand{\bm}[1]{\boldsymbol{#1}}
\newcommand{\vect}[1]{\bm{#1}}
\newtheorem{remark}{Remark}[section]
\newcommand{\reffig}[1]{\text{Fig.}\ref{#1}}
\newcommand{\unit}[1]{\mathrm{#1}}

\title[3D Finite Element Modeling with Impact Ionization]
{3D Finite Element Modeling of Current Densities in 
Semiconductor Transport with Impact Ionization}

\address{$^1$ Modeling and Simulations group
         Micron Technology 
         via Olivetti 2, Agrate Brianza, Italy \\
				 $^2$ IBM Italia, 
Via Circonvallazione Idroscalo, 20090, Segrate (MI), Italy \\
         $^3$ STMicroeletronics, via Olivetti 2, 20864, 
Agrate Brianza (MB), Italy \\ 
				 $^{4}$ Dipartimento di Matematica,
         Politecnico di Milano 
				 Piazza L. da Vinci 32, 20133 Milano, Italy}

\author{A. Mauri$^{1}$ \and A. Bortolossi$^2$ \and
G. Novielli$^3$ \and R. Sacco$^{4}$}
\email{aureliogianc@micron.com} 
\email{andrea.bortolossi@it.ibm.com}
\email{giovanni.novielli@st.com}
\email{riccardo.sacco@polimi.it}

\begin{document}

\date{\today}

\begin{abstract}
In this article we propose two novel 3D finite element models, 
denoted method A and B, for electron and hole Drift-Diffusion (DD)
current densities.
Method A is based on a primal-mixed formulation of the 
DD model as a function of the quasi-Fermi potential gradient, while 
method B is a modification of the standard DD formula
based on the introduction of an artificial diffusion matrix.
Both methods are genuine 3D extensions of the classic
1D Scharfetter-Gummel difference formula. 
The proposed methods are compared in the 3D simulation
of a $p$-$n$ junction diode and of a $p$-MOS transistor
in the on-state regime. 
Results show that method A provides the best performance in terms of physical accuracy and numerical stability. 
Method A is then used in the 3D simulation of a $n$-MOS transistor in the off-state regime including the impact ionization 
generation mechanism. Results demonstrate that the model is able to accurately compute the I-V characteristic of the device until drain-to-bulk junction breakdown.
\end{abstract}

\maketitle

{\bf Keywords:}
Drift-Diffusion model; Finite Element method;
Numerical Simulation; Impact Ionization;
PACS:85.30; MSC 2000: 35J; 65N30.

\section{Introduction and Motivation}\label{sec:intro}

Semiconductor technology is undergoing a continuously increasing
advancement in the aggressive scaling of device length~\cite{ITRS2013}.
In this scenario, three-dimensional (3D) device modeling and numerical simulation techniques
play a critical role in 
the prediction of
the electrical performance of the system under investigation.
The result of our approach to 3D modeling and simulation of
semiconductor device applications (see~\cite{NovielliIEDM2013} and~\cite{Mauri2014}) is the 
FEMOS project (Finite Element Modeling Oriented Simulator), 
a general-purpose modular numerical code based on the use of the 
Galerkin Finite Element Method (GFEM) implemented in a fully 3D framework 
through shared libraries using an objected-oriented programming 
language (C++).
In the present work we employ the FEMOS computational platform in 
the study of the \textit{Drift Diffusion model} 
(DD)~\cite{selberherr:SimSem,Jerome:AnalyCharTran} 
and focus on the issue of endowing the simulation tool
of a consistent, stable and accurate procedure for the approximation
of electron and hole current densities in the device.
This is of utmost importance in: 
i) visualization and post-processing; 
ii) evaluation of conduction currents at device terminals; 
and iii) inclusion in the DD model of generation phenomena due to
Impact Ionization (II). Here, our attention is 
devoted to iii), because of the critical role of II in the 
convergence and numerical stability of the iterative algorithm 
used to solve the DD system (see~\cite{markowich1986stationary}, 
Chapt. 3 and~\cite{Micheletti1995}), although the methods we propose for
the treatment of iii) can also be profitably employed for i) and ii).

To allow a consistent treatment of the generation term due to II within the
FE procedure, we propose two novel discrete 
models for electron and hole current densities over the computational grid.
The two methods provide a constant approximation of the current density inside each mesh element and for this reason they can be easily 
implemented in any simulation environment
not necessarily employing the GFEM but utilizing, instead, 
the more standard Box Integration Method 
(see~\cite{BankRoseFichtnersiam1983,BankRoseFichtnerieee1983} and~\cite{selberherr:SimSem}, Chapt. 6).
The first scheme (method A) is based on the use of a primal-mixed 
formulation of the DD model written as a function of the quasi-Fermi 
potential gradient (see~\cite{Roberts:MixedHybrid} and~\cite{Sacco2014ABB}). 
The second scheme (method B) is a modification
of the standard DD formula based on the introduction of an artificial diffusion 
matrix (see also~\cite{Bank:Upwinding}).
Both methods 
are genuine 3D extensions of the classic
Scharfetter-Gummel (SG) formula for the computation of the current density 
over a 1D element~\cite{Gummel:SignAnalys}. 

The proposed finite element models are validated 
in the numerical study of 3D device structures
($p$-$n$ junction diode and MOS transistors) under on-state
and off-state regimes. 
Results clearly indicate that method A provides 
the best performance in terms of physical accuracy and numerical stability, and 
demonstrate the ability of 
the simulation model in accurately computing the I-V characteristic of the device until the onset of drain-to-bulk junction breakdown.

The outline of the article is as follows. 
In Sect.~\ref{sec:models_geo_discretization} 
we briefly review the DD model 
and its finite element discretization. 
In Sect.~\ref{sec:fem} we illustrate the novel methods 
proposed to calculate the current densities in the device.
In Sect.~\ref{sec:numexp} we carry out their extensive validation by comparing the obtained results with a reference simulation suite in the study of 3D $p$-$n$ and MOS structures including the II generation term.
Concluding remarks and perspectives are addressed in 
Sect.~\ref{sec:conclusion}.

\section{Transport Model and Discretization}
\label{sec:models_geo_discretization}
The Drift-Diffusion (DD) equations for 
electron and hole current densities $\vect{J}_n$ and $\vect{J}_p$ 
in a semiconductor device are:
\begin{subequations}\label{eq:DDcurrents}
\begin{align}
& \vect{J}_n = q \mu_n n \vect{E} + q D_n \nabla n  & \label{eq:Jn_vector} \\
& \vect{J}_p =  q \mu_p p \vect{E} - qD_p \nabla p & \label{eq:Jp_vector} 
\end{align}
\end{subequations}
where $n$ and $p$ are the densities of electrons and holes, 
$\vect{E}=-\nabla \varphi$ is the electric field, $\varphi$ being the
electric potential, while $\mu_n$ and $\mu_p$ 
are the electron and hole mobilities and $D_n$ and $D_p$
are the electron and hole diffusivities, related to the mobilities
through Einstein's relation. We refer 
to~\cite{selberherr:SimSem,markowich1986stationary,Jerome:AnalyCharTran}
for an extensive analysis and discussion of the DD model.
To iteratively solve the differential equation system constituted by
the Poisson equation for $\varphi$ and 
the continuity equations for $n$ and $p$
we adopt the Gummel decoupled solution map with the lagging 
technique for the treatment of recombination and generation mechanisms.
We refer to~\cite{GummelMap} and~\cite{Jerome:AnalyCharTran} for
a description and thourough analysis of the solution map.

The device under investigation is geometrically 
represented by the polyhedral 
domain $\Omega \subset \R^3$, given by the union of two open disjoint subsets, a doped silicon part and an oxide part 
assumed to be a perfect insulator, separated by an interfacial surface.
For the numerical simulation of the DD model, 
the domain $\Omega$
is divided in a discrete partition $\mathcal{T}_h$ made by elements $K$,
each element being a tetrahedron of volume $|K|$, such that 
$\overline{\Omega} =  \bigcup_{K \in \mathcal{T}_h} \overline{K}$. 
Then, each differential problem in the Gummel decoupled algorithm 
is written in weak form 
(see~\cite{quarteroni:NumApprox}, Chapt. 5) and discretized
using the displacement-based 
Galerkin Finite Element Method (GFEM) with piecewise linear conforming 
elements for potential and carrier densities 
(see~\cite{quarteroni:NumApprox}, Chapt. 6). 
To avoid numerical
instabilities due to possible dominance of the drift term,
the variant of the GFEM denoted Edge Averaged Finite Element (EAFE) method 
proposed and analyzed in~\cite{Zikatanov:EAFE1,Zikatanov:EAFE2} is used 
in the solution of the linearized continuity equations for electrons
and holes.

\section{Finite Element Models for the Current Density}\label{sec:fem}
 
The construction of a stable and accurate 
approximate current density field 
in a primal-based FE formulation is not a trivial task because
of possible numerical problems arising from differentiation and cancellation
in the DD transport relations~\eqref{eq:Jn_vector}-\eqref{eq:Jp_vector}.
In this section we introduce two novel finite element methodologies 
for current density discretization. The proposed approaches 
have a cheap computational cost compared to other 
formulations (such as the dual-mixed FEM, see~\cite{Brezzi1989}), 
are completely compatible for use in the classical Box Integration 
Method (BIM~\cite{BankRoseFichtnersiam1983,BankRoseFichtnerieee1983})
and are genuine 3D extensions of the classic Scharfetter-Gummel
(SG) 1D difference formula~\cite{Gummel:SignAnalys}.
In the remainder of the 
article, for a given element $K$ in $\mathcal{T}_h$ we denote the volume of 
$K$ by ${\rm vol}(K)$; moreover, with the subscript $K$ we refer to a quantity defined in the interior of $K$ while the subscript $h$ refers to a quantity defined at the vertices of $K$. 
For a given function $f : K \rightarrow \R$ we define
$\langle f \rangle_K:=\int_K f \, dK /{\rm vol}(K)$ the mean
integral value of $f$ over $K$.
We also assume that carrier mobility (and the associated diffusivity through
Einstein's relation) is constant in $K$, while the electric potential is linear
(so that the associated electric field is a constant vector in $K$).

\subsection{The DD Formula}\label{sec:method_DD}

The simplest FE model for the current densities (method DDFE) is 
obtained by replacing into the transport equations~\eqref{eq:Jn_vector}-\eqref{eq:Jp_vector} the functions
$n$, $p$ and $\varphi$ with their corresponding FE approximations
$n_h$, $p_h$ and $\varphi_h$, and then by computing the integral average
of the resulting expressions over the element $K$.
This yields
\begin{subequations}\label{eq:method_A}
\begin{align}
& \vect{J}_{n,K} = q \mu_n \langle n \rangle_{K} \vect{E}_K + q D_n \nabla n_{h,K} 
\label{eq: Jn DD discrete}\\ 
& \vect{J}_{p,K} = q \mu_p \langle p \rangle_{K} \vect{E}_K - 
q D_p \nabla p_{h,K}.
\label{eq: Jp DD discrete} 
\end{align}
\end{subequations}
It is immediate to see that the discrete current densities~\eqref{eq: Jn DD discrete} 
and~\eqref{eq: Jp DD discrete} automatically 
reproduce the limiting cases of pure diffusive flow ($\vect{E}_K = \vect{0}$)
and pure ohmic flow ($n_{h,K} = {\rm constant}$ and 
$p_{h,K} = {\rm constant}$). In the case of thermal equilibrium
($\vect{J}_{n,K} = \vect{J}_{p,K} = \vect{0}$), we can anticipate 
computational difficulties with the use of~\eqref{eq: Jn DD discrete} 
and~\eqref{eq: Jp DD discrete} because of exact cancellation
of drift and diffusive current contributions. Thus, by a continuity argument,
we also see that method DDFE does not seem appropriate in the numerical treatment of the subthreshold current regime, where currents are not exacty equal to zero but are very small. The two following formulations are
designed to overcome this limitation.

\subsection{Method A}\label{sec:method_A}
To describe method A we consider the case of electron continuity
equation because a completely similar treatment holds for the
hole continuity equation. We start from the Maxwell-Boltzmann (MB) 
statistics for electrons 
$n = n_i \exp \left( (\varphi_n-\varphi)/V_{th} \right)$
where $n_i$ is the intrinsic concentration in the semiconductor 
material, $\varphi_n$ is the quasi-Fermi potential for electrons and
$V_{th}$ is the thermal voltage. Replacing the MB relation
into~\eqref{eq:Jn_vector} yields the equivalent (nonlinear) 
gradient form of the electron current density 
\begin{subequations}\label{eq:method_B}
\begin{align}
& \vect{J}_n = - q \mu_n n \nabla \varphi_n = 
- q \mu_n n_i \exp \left( (\varphi_n-\varphi)/V_{th} \right)
\nabla \varphi_n. & 
\label{eq:J_n_ohmic}
\end{align}

To construct the finite element model for $\vect{J}_n$ as in~\eqref{eq:J_n_ohmic}, we use the primal-mixed (PM) FEM introduced
and analyzed in~\cite{Roberts:MixedHybrid} and recently extended to the case of advective-diffusive operators in~\cite{Sacco2014ABB}.
In the PM FEM of lowest order, 
the approximate current density is constant over each $K \in \mathcal{T}_h$
while the approximate quasi-Fermi potential is piecewise linear and 
continuous over $\mathcal{T}_h$. 
Let us introduce the finite element spaces of piecewise constant and
piecewise linear continuous functions over $\mathcal{T}_h$:
\begin{align}
& Q_h=\left\{ w \in L^2(\Omega) : w|_{K}\in \P_0(K) \forall K \in \mathcal{T}_h\right\} & \label{eq:Qh} \\
& V_h=\left\{ u \in C^0(\overline{\Omega}) : u|_{K} \in 
\P_1(K) \forall K \in \mathcal{T}_h \right\} & \label{eq:Vh}
\end{align}
and the electrical conductivity 
$\sigma_{n,h}:= q \mu_n n_i \exp \left( (\varphi_n-\varphi)/V_{th} \right)$.
Then, the
PM-FE approximation of~\eqref{eq:J_n_ohmic} reads: 
find $\vect{J}_{n,h} \in \left[Q_h\right]^3$ such that
\begin{align}
& 
\int_\Omega 
\sigma_{n,h}^{-1} \vect{J}_{n,h} \cdot \vect{q}_h \, d\Omega
+ \int_\Omega \nabla \varphi_{n,h} \cdot \vect{q}_h \, d\Omega = 0 
& \forall \vect{q}_h \in \left[Q_h\right]^3 \label{eq:PMFE_Jh} 
\end{align}
where $\vect{J}_{n,h} \in \left[Q_h\right]^3$, 
$\varphi_{n,h} \in V_h$ and $\varphi_{h} \in V_h$
are the finite element discrete analogues of $\vect{J}_{n}$, 
$\varphi_n$ and $\varphi$. 
Since functions in $\left[Q_h\right]^3$ are discontinuous, 
relation~\eqref{eq:PMFE_Jh} amounts to
\begin{align}
& 
\int_K 
\sigma_{n,h}^{-1} \vect{J}_{n,h} \cdot \vect{q}_h \, dK
+ \int_K \nabla \varphi_{n,h} \cdot \vect{q}_h \, dK = 0 
& \forall \vect{q}_h \in \left[\P_0(K)\right]^3. \label{eq:PMFE_Jh_K} 
\end{align}
Using in~\eqref{eq:PMFE_Jh_K} the standard basis functions 
for $\left[\P_0(K)\right]^3$ 
we obtain
\begin{align}
&
[\mathbf{J_{n,h}}]_i = 
- \mathcal{H}_K \left( \sigma_{n,h} \right)  
\dfrac{\partial \varphi_{n,h}}{\partial x_i} \quad i = 1,2,3, 
\quad \forall K \in \mathcal{T}_h
& \label{eq: first formula for J}
\end{align}
where $\mathcal{H}_K \left( \sigma_{n,h} \right) := \left( \langle \sigma_{n,h}^{-1} \rangle_K \right)^{-1}$ 
is the harmonic average of $\sigma_{n,h}$ over the element $K$.

To numerically compute in a simple and accurate manner the harmonic
average of the electrical conductivity, we use 
the following quadrature 
rule (see also~\cite{Brezzicmame1989,BrezziLectNotes1991})
\begin{align}
& \left(\dfrac{\int_K \sigma^{-1}_{n,h} \, dK}{{\rm vol}(K)} 
\right)^{-1} \simeq \left(\dfrac{\int_{e~\ast} 
\sigma^{-1}_{n,h} \, de}{|e^\ast|} \right)^{-1}
& \label{eq: approximation from 3D to edge}
\end{align}
$e^*$ being the edge of $\partial K$ where 
the maximum drop of $\sigma_{n,h}$ occurs and $|e^\ast|$ its
euclidean length.
To identify the edge $e^*$ we introduce the linear dimensionless
potential $\Phi := \left(\varphi_h - 	\varphi_{n,h}\right)/V_{th}$
and we determine the two vertices: $\mathbf{x}_m$ s.t. $\Phi(\mathbf{x}_m)=\Phi_m := \min_K(\Phi)$ and $\mathbf{x}_M$ s.t. $\Phi(\mathbf{x}_M)=\Phi_M:=\max_K(\Phi)$. Then, we define 
$e^\ast$ to be the edge of $\partial K$ which connects 
$\mathbf{x}_m$ and $\mathbf{x}_M$.
The evaluation of the approximate 1D integral 
in~\eqref{eq: approximation from 3D to edge}
by choosing the orientation from $\mathbf{x}_m$ and $\mathbf{x}_M$ yields
\begin{align}
& \int_{K} \sigma_{n,h}^{-1} \, dK \simeq  
q \mu_n n_i \exp(\Phi_m) \mathcal{B}(\Phi_m-\Phi_M) & 
\label{eq: finally approximation 3D to 1D}
\end{align}
while choosing the orientation from $\mathbf{x}_M$ and $\mathbf{x}_m$
the result is
\begin{align}
& \int_{K} \sigma_{n,h}^{-1} \, dK \simeq  
q \mu_n n_i \exp(\Phi_M) \mathcal{B}(\Phi_M-\Phi_m), & 
\label{eq: approssimazione sm}
\end{align}
where $\mathcal{B}(Z): = Z/(\exp(Z) -1)$ 
is the inverse of the Bernoulli function, such that
$\mathcal{B}(0)=1$.
Using MB statistics into the previous two
relations, the property $e^Z \mathcal{B}(Z) = \mathcal{B}(-Z)$, and
combining (\ref{eq: finally approximation 3D to 1D}) and (\ref{eq: approssimazione sm}), we find
\begin{align}
& \mathbf{J}_{n,K} = -  q \mu_n  
\left[ \dfrac{ n_m \mathcal{B}(-\Delta \Phi_{max})  + 
n_{M}\mathcal{B}(\Delta \Phi_{max})}{2} 
\right]\nabla \varphi_{n,h} & \label{eq: second formula for J}
\end{align}
where $n_m: = n_i e^{\Phi_m}$ and $n_M :=n_i e^{\Phi_M}$ while 
$\Delta \Phi_{\max} := \Phi_M - \Phi_m$. 

\begin{remark}[Method A and the 1D SG scheme]\label{rem:method_B_SG}
The approximate electron current density~\eqref{eq: second formula for J} 
can be regarded as a 
consistent extension to the 3D case of the classic 
1D SG formula, provided to replace over the tethrahedron $K$ the 
1D gradient of the electron quasi Fermi potential $\varphi_n$ 
with the constant 3D gradient of the 
approximation $\varphi_{n,h}|_K$, and the 
electrical conductivity $\sigma_n$ with the constant 
electrical conductivity $\sigma_{n,h}\big|_K$ 
defined (as in the SG formula) as the 1D harmonic average of 
$\sigma_n$ over the interval 
$[\mathbf{x}_m, \, \mathbf{x}_M]$. 
\end{remark}
Following the same procedure also for the hole current density we have
\begin{align}
\mathbf{J}_{p,K} = -  q \mu_p  
\left[ \dfrac{ p_m \mathcal{B}(\Delta \Phi_{max}) + 
p_M \mathcal{B}(-\Delta \Phi_{max})}{2} \right]\nabla \varphi_{p,h}.
& \label{eq: first formula for Jp}
\end{align}
\end{subequations}

\subsection{Method B}\label{sec:method_B}
In the previous section the discrete form of the 
current density is constructed by starting from the equivalent 
"ohmic" representation in terms of the quasi Fermi potential, and then 
by performing a suitable approximation of the 
electrical conductivity over the finite element $K$.
Here, we continue along the same line of thought, but starting from
the classic DD format~\eqref{eq:Jn_vector}-~\eqref{eq:Jp_vector}, with
the intent of using the method of Streamline Upwind 
artificial diffusion proposed 
in~\cite{Brooks1982} for the advective-diffusive model to stabilize 
the computation in the presence of a high electric field.

\subsubsection{The 1D SG Method as an Artificial Diffusion Scheme}\label{sec:method_1d_SG}
In the 1D setting the artificial diffusion 
technique consists of replacing the electron
diffusion coefficient $D_{n}$ with the modified quantity
\begin{subequations}\label{eq:method_C}
\begin{align}
& D_{n,h} = 
D_n + D_n \Phi \left(\P{}e|_K\right) & \label{eq: perturbed diffusion}
\end{align}
where $\Phi$ is a suitable nonnegative stabilization function of 
the 1D local P\`eclet number 
$\P{}e|_K = \left(h |\partial_x \varphi_h|\right)/(2V_{th}) = 
|\Delta \varphi|/(2 V_{th})$,  
$h$ and $\Delta \varphi$ being the length of the 1D interval and the
potential drop over the interval, respectively.
The local P\`eclet number gives a measure of how much the drift term
dominates over the diffusion term in the transport mechanism.
If $\P{}e|_K > 1$ the problem is locally drift (advection)-dominated 
and in such a case we need introduce an extra amount of diffusion
in~\eqref{eq: perturbed diffusion} (given by $D_n \Phi(\P{}e|_K)$) 
to prevent the occurrence of unphysical spurious oscillations in the
computed solution, which might even lead to 
a negative electron concentration. If $\P{}e|_K < 1$ 
the problem is locally diffusion-dominated and there is no need
of adding an extra diffusion, so that the standard GFEM is enough
for obtaining an accurate and numerically stable solution.
Based on the last observation, 
the function $\Phi$ has to satisfy the property of consistency
\begin{align}
\label{eq: consistenza}
\lim_{\P{}e|_K \to 0} \Phi(\P{}e|_K) = 0 \quad\quad \forall K\in\mathcal{T}_h.
\end{align}
The 1D approximation of the electron current density to be used in a 
stabilized GFEM is thus given by the following relation
\begin{align}
& J_{n,h}(n_h) \big|_K = q \mu_n \langle n \rangle_K 
E_{h,K} + q D_n
\left( 1 + \Phi(\P{}e|_K)) \right)
\partial_x n_h & \forall K\in\mathcal{T}_h. \label{eq:Jn_1D_art_diff}
\end{align}
To design in a physically sound and consistent manner the 
optimal stabilization function $\Phi$, we pretend the 
modified method to exactly satisfy some limiting cases 
that often occur in practical important electronic applications.
Using a 3D notation, for sake of generality and because this will
be used in later extension, the considered cases are:
\begin{itemize}
\item[C1.] \textbf{Constant carrier concentrations} (only drift contribution):
$\vect{J}_n  = q\mu_n n \vect{E}.$
\item[C2.] \textbf{Constant potential} ($\vect{E}=0$, only diffusive contribution): $\vect{J}_n  = qD_n\nabla n.$
\item[C3.] \textbf{Constant quasi Fermi potential} (no current flow): 
$\vect{J}_n = -q \mu_n n \nabla \varphi_n = \vect{0}$.
\end{itemize}
Notice that case C3. implies that 
\begin{align}
& n=C e^{\varphi / V_{th}} & \label{eq:n_exponential}
\end{align}
where $C$ is an arbitrary constant such that
$C = \exp(-\overline{\varphi}_n/V_{th})$, 
$\overline\varphi_n$ being a given constant value.
Thanks to assumption~\eqref{eq: consistenza} the stabilized current 
\eqref{eq:Jn_1D_art_diff} automatically satisfies cases C1 and C2. 
Case C3 is recovered by imposing (in the 1D setting)
\begin{equation}
\label{eq: condizione JnK nulla}
J_{n,h}(\Pi_1^K(Ce^{\varphi / V_{th}})) = 0
\end{equation}
where $\Pi_1^K$ is the $\P_1$-interpolant over the element $K$.
Using \eqref{eq: condizione JnK nulla} in 
\eqref{eq:Jn_1D_art_diff}, noting that $E_K = -\partial_x 
\varphi_h$ and using Einstein's relation, 
we obtain 
the following relation for the stabilization function
\begin{equation}
\label{eq: SG stab func alternative}
\Phi(\P{}e|_K)  = 
\dfrac{\langle n \rangle_K}{V_{th}} 
\dfrac{\partial_x \varphi_h}{\partial_x \Pi_1^K(e^{\varphi / V_{th}}) } -1. 
\end{equation}
Enforcing relation~\eqref{eq:n_exponential} at the two vertices $x_i$ 
of the interval, $i=1,2$, yields
\begin{align}
\Phi(\P{}e|_K) & = \sigma \P{}e|_K \dfrac{e^{\varphi_1/V_{th}}+e^{\varphi_2/V_{th}}}
{e^{\varphi_2/V_{th}}-e^{\varphi_1/V_{th}}} -1 
= \sigma \P{}e|_K 
\dfrac{e^{2 \sigma \P{}e|_K}+1}{e^{2\sigma \P{}e|_K}-1} -1 & 
\label{eq:Phi_stab}
\end{align}
where $n_1$ and $n_2$ are the two nodal values of $n_h$ while
$\sigma: = {\rm sign}(\Delta \varphi)$.
Setting for brevity $X := 2 \sigma \P{}e|_K$, 
relation~\eqref{eq:Phi_stab} becomes
\begin{align*}
\Phi(X) & = \dfrac{X}{2} \left( \dfrac{e^{X}}{e^{X}-1} + \dfrac{1}{e^{X}-1} \right) -1 = \dfrac{1}{2} \left( \mathcal{B}(-X) + 
\mathcal{B}(X) \right) -1 \\
& = \dfrac{1}{2} \left( X + \mathcal{B}(X) + \mathcal{B}(X) \right) -1 
=  \mathcal{B}(X) +\dfrac{X}{2}  -1,
\end{align*}
and replacing the definition of $X$ we 
obtain for both $\Delta \varphi > 0$ and $\Delta \varphi < 0$
\begin{equation}
\label{eq: SG stabilization function}
\Phi(\P{}e|_K) = \mathcal{B}(2\P{}e|_K) + \P{}e|_K -1
\end{equation}
which, upon substitution into~\eqref{eq:Jn_1D_art_diff},  
recovers the well known 1D Scharfetter-Gummel scheme 
originally also proposed as "exponential fitting"
by Allen and Southwell in~\cite{Allen:SGformula}.

\subsubsection{The 3D SG Artificial Diffusion Method}
\label{sec:method_3d_SG}
In a 3D framework a straightforward extension of the 1D Scharfetter-Gummel stabilization~\eqref{eq:Phi_stab}-~\eqref{eq: SG stabilization function} 
can be obtained by introducing a $3\times 3$ diagonal stabilizing tensor  $\underline{\underline{\vect{\Phi}}}^K$ defined on each element $K
\in \mathcal{T}_h$ as follows
\begin{equation}
\label{eq: SG stabilized}
\underline{\underline{\vect{\Phi}}}^K_{\, ii}  = \dfrac{\langle \Pi_1^K(e^{(\varphi_h-\varphi_{M})/V_{th}}) \rangle_K \partial_{x_i}\varphi}{\partial_{x_i}\Pi_1^K
(e^{(\varphi_h-\varphi_{M})/V_{th}})V_{th}} -1 \quad\quad i = 1,2,3.
\end{equation}
where we have used, as a reference value for the potential reference
to avoid overflow exceptions in the machine 
evaluation of~\eqref{eq: SG stabilized}, 
the maximum $\varphi_{M}$ of $\varphi_h$ in $K$.

The 3D approximate electron current density is then
\begin{equation}
\label{eq: Jn 3D upwin SG}
\vect{J}_{n,K} = q \mu_n \langle n \rangle_K \vect{E}_K 
+ qD_n(\underline{\underline{\mathcal{I}}}+\underline{\underline{\vect{\Phi}}}_K) \nabla n_h
\end{equation}
where $\underline{\underline{\mathcal{I}}}$ is the $3\times 3$ identity tensor.
In a completely similar manner we have the relation for 
hole current density
\begin{equation}
\label{eq:  Jp 3D upwin SG}
\vect{J}_{p,K} = q \mu_p \langle p \rangle_K \vect{E}_K - qD_p(\underline{\underline{\mathcal{I}}}+
\underline{\underline{\vect{\Phi}}}_K) \nabla p_h
\end{equation}
where
\begin{equation}
\underline{\underline{\vect{\Phi}}}^K_{\, ii}  = \dfrac{\langle \Pi_1^K(e^{(\varphi_{m}-\varphi_h)/V_{th}})\rangle_K \partial_{x_i}\varphi}{\partial_{x_i}
\Pi_1^K(e^{(\varphi_{m}-\varphi_h)/V_{th}})V_{th}} -1 \quad\quad i = 1,2,3,
\end{equation}
$\varphi_{m}$ being the minimum of $\varphi$ over $K$.
It is immediate to check that the two proposed novel 
approximations~\eqref{eq: Jn 3D upwin SG} and~\eqref{eq:  Jp 3D upwin SG} satisfy all cases C1, C2 and C3 in Sect.~\ref{sec:method_1d_SG}.
\end{subequations}

\section{Numerical Experiments}
\label{sec:numexp}
In order to compare the performance of the methods of Sect.~\ref{sec:fem} we
use the FEMOS computational platform for the numerical simulation of 
three 3D devices: 1) a $p$-$n$ diode; 2) a $n$-MOS device; and 3) a 
$p$-MOS device in on and off working conditions. 
While the diode and the $n$-MOS devices are coming out from ideal structures, the $p$-MOS is the result of a realistic 2D-process simulation accounting 
for non ideal doping profiles. This last test structure is a severe benchmark
to highlight the accuracy and the stability of the different density current calculation methods.   

\subsection{Diode}\label{diode}
\begin{figure}[!h]
\centering
\subfloat[][\emph{Mesh}\label{fig: pn mesh}]
{\includegraphics[width = 0.45\textwidth , height=4.5cm]{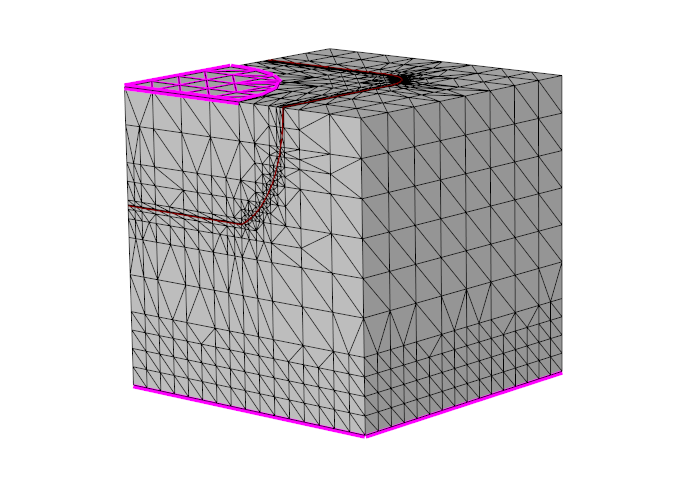}}
\hspace{0.05\textwidth}
\subfloat[][\emph{2D cut of Mesh for doping profile visualization}\label{fig: pn doping}]
{\includegraphics[width = 0.45\textwidth , height=4.5cm]{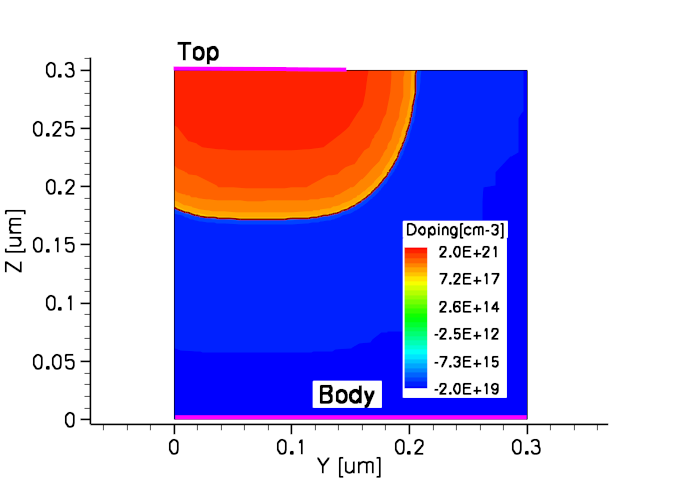}}
\caption{Diode structure: Left: doping. Right: mesh.}
\label{fig: pn struct}
\end{figure}

\begin{figure}[!h]
\centering
\subfloat[][\emph{Method DDFE - $\vect{J}_n$}\label{fig: Jn DD}]
{\includegraphics[width = 0.5\textwidth, height=4.5cm]{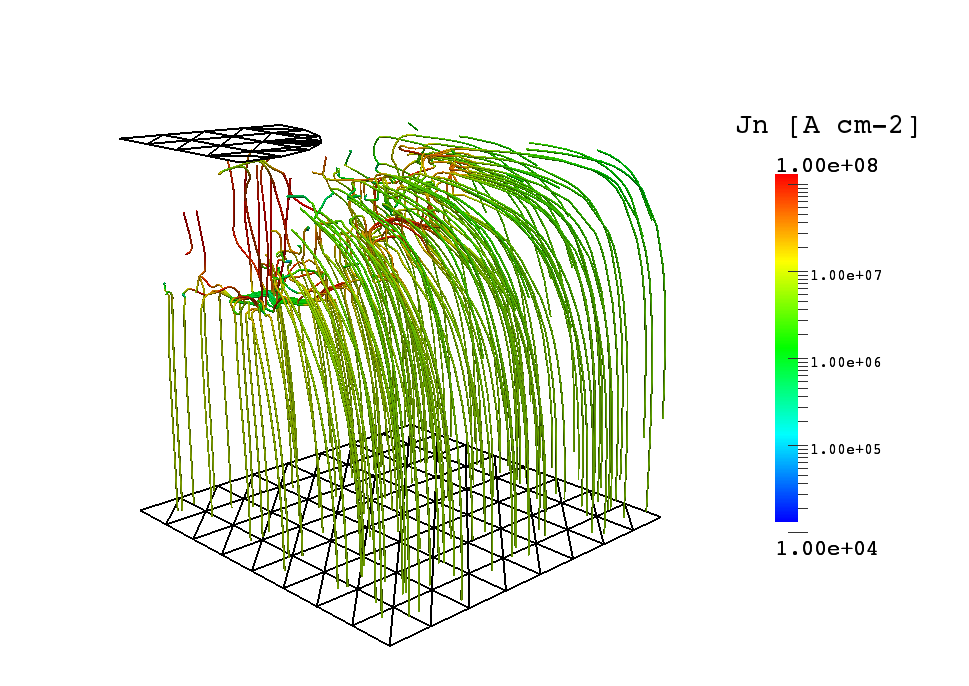}}
\subfloat[][\emph{Method DDFE - $\vect{J}_p$}\label{fig: Jp DD}]
{\includegraphics[width = 0.5\textwidth , height=4.5cm]{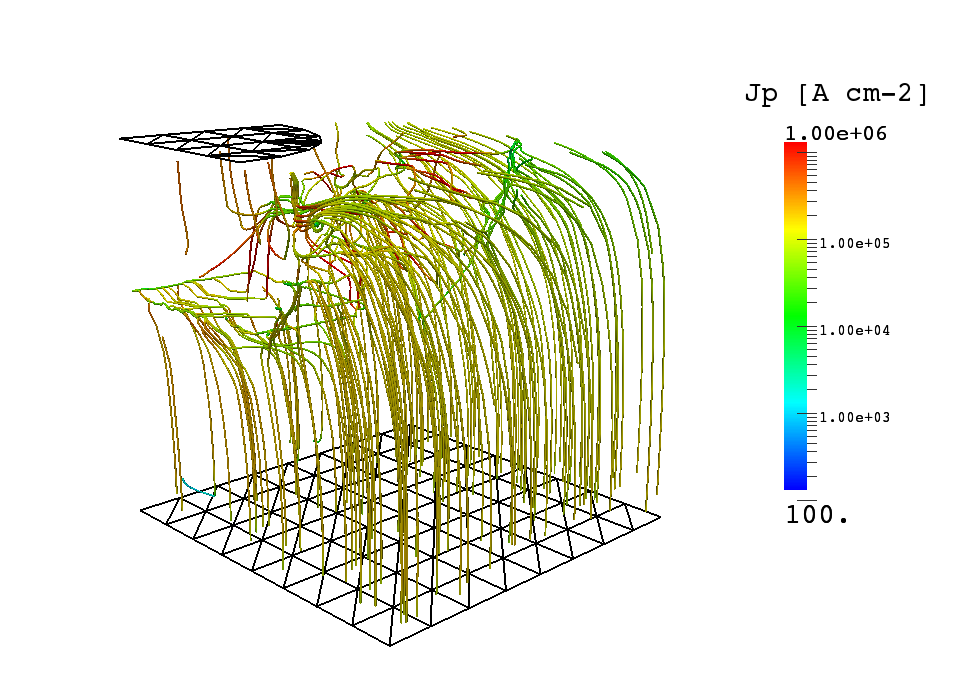}} \\
\subfloat[][\emph{Method A - $\vect{J}_n$}\label{fig: Jn SG}]
{\includegraphics[width = 0.5\textwidth , height=4.5cm]{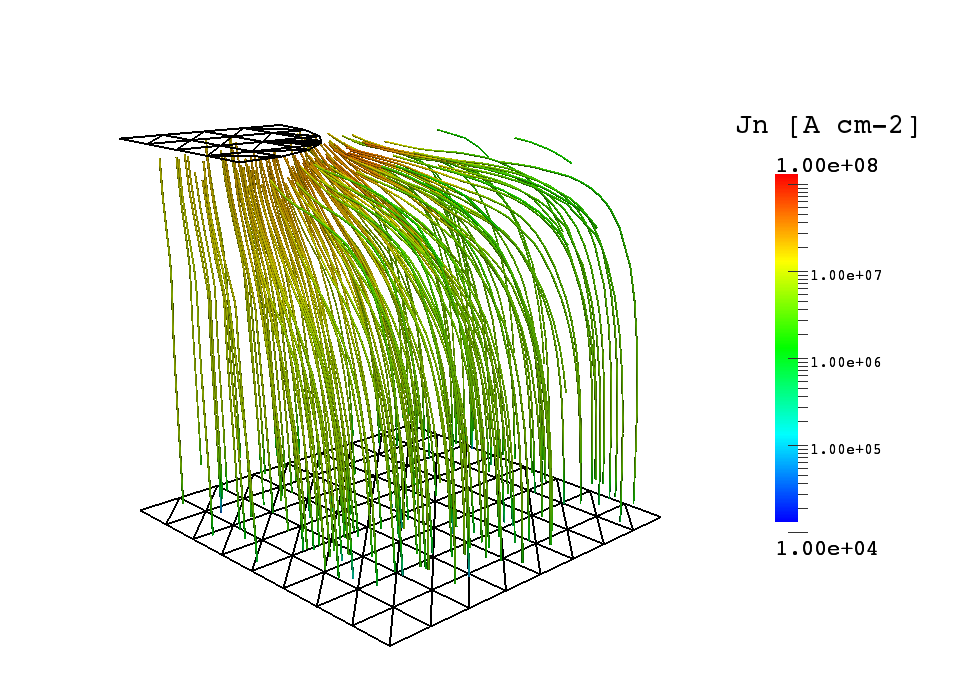}}
\subfloat[][\emph{Method A - $\vect{J}_p$}\label{fig: Jp SG}]
{\includegraphics[width = 0.5\textwidth , height=4.5cm]{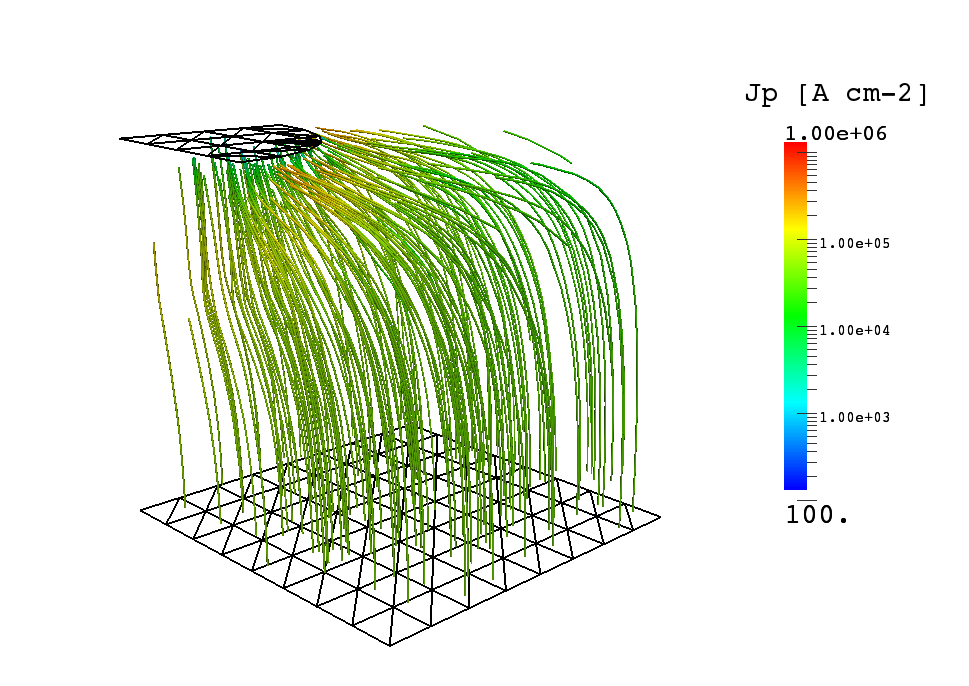}} \\
\subfloat[][\emph{Method B - $\vect{J}_n$}\label{fig: Jn UPWIND}]
{\includegraphics[width = 0.5\textwidth , height=4.5cm]{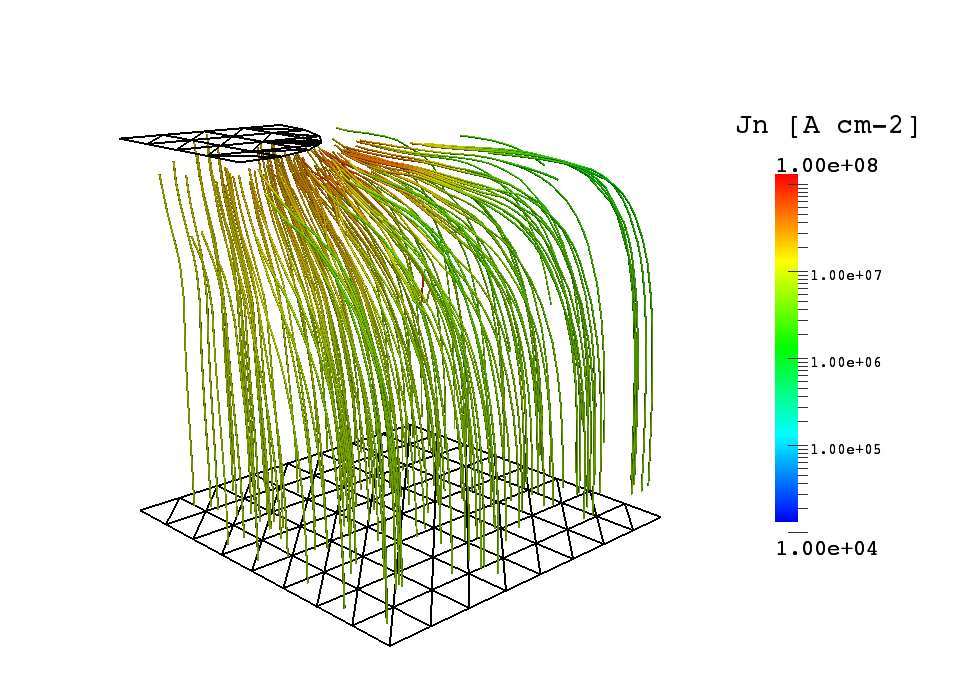}}
\subfloat[][\emph{Method B - $\vect{J}_p$}\label{fig: Jp UPWIND}]
{\includegraphics[width = 0.5\textwidth , height=4.5cm]{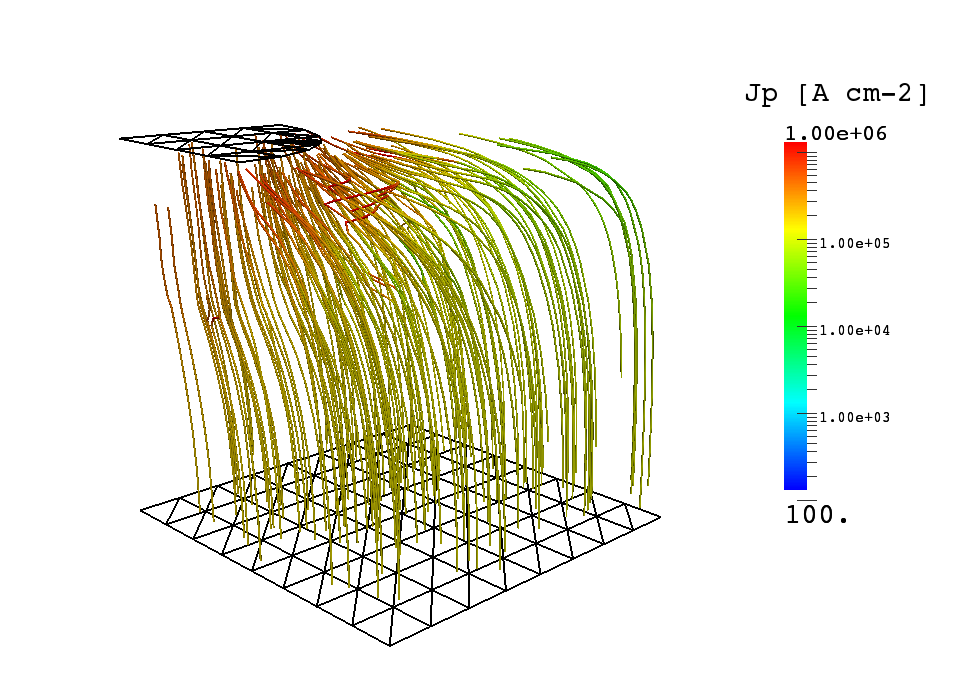}}
\caption{Electron (left column) and hole (right column) 
current densities of Diode test case at $V_{body}=0.8[V]$. 
Top: Method DDFE. Middle: Method A. Bottom: Method B.}
\label{fig: pn current density 0,8V}
\end{figure}

The first structure is a semiconductor region with a $p$-$n$ junction whose dimensions are $\Omega=(0,0.3)^3 {\mu m}^3$. 
A Gaussian implantation of donors with a peak of 
$10^{19}[cm{-3}]$ and a depth of $0.15[\mu m]$ 
\reffig{fig: pn doping} is made over an $p$-type 
region with a constant acceptor profile of $10^{18}[cm^{-3}]$ magnitude. 
Two contacts are defined for each of the doped regions: for the 
$n$-type part a rounded-shape contact is used (Top), while the $p$-type 
part is contacted at the bottom face (see \reffig{fig: pn mesh}) (Body).
Contacts are considered as a pure Ohmic-type with appropriate 
Dirichlet boundary conditions: the Top is maintained at ground while 
the Body is ramped at $0.8[V]$. 
In \reffig{fig: pn current density 0,8V} the results of the calculation of the current density flux for electrons and holes in the semiconductor bulk are represented through streamlines connecting the Body with the Top contact.  
As expected the current calculation obtained with method DDFE (eq.~\eqref{eq:method_A}) is affected by a critical behavior, in particular close to and inside the $n$-junction as shown in~\reffig{fig: Jn DD} and~\reffig{fig: Jp DD} where instability has to be ascribed to numerical cancellation of the drift and diffusion contributions. 
Results get definitely better by employing method B using
Eqns.~(\ref{eq: Jn 3D upwin SG}) and (\ref{eq:  Jp 3D upwin SG}) where the improvement can be appreciated in \reffig{fig: Jn UPWIND} and \reffig{fig: Jp UPWIND}. However, 
a careful inspection of the hole current density reveals that some numerical instability is still evident inside the $n$-junctions. 
The extension of the 1D SGscheme to 3D provided by method A in Eq.(\ref{eq: first formula for J}) and (\ref{eq: first formula for Jp}) results in the streamlines presented in \reffig{fig: Jn SG} and \reffig{fig: Jp SG}: 
no spurious instability can be observed anymore and our calculations are in excellent agreement with the results of a commercial code (not shown here).

\subsection{$p$-MOS structure}

\begin{figure}[!h]
\centering
\subfloat[][\emph{Doping}\label{fig: MOS doping}]
{\includegraphics[width = 0.45\textwidth , height=4.5cm]{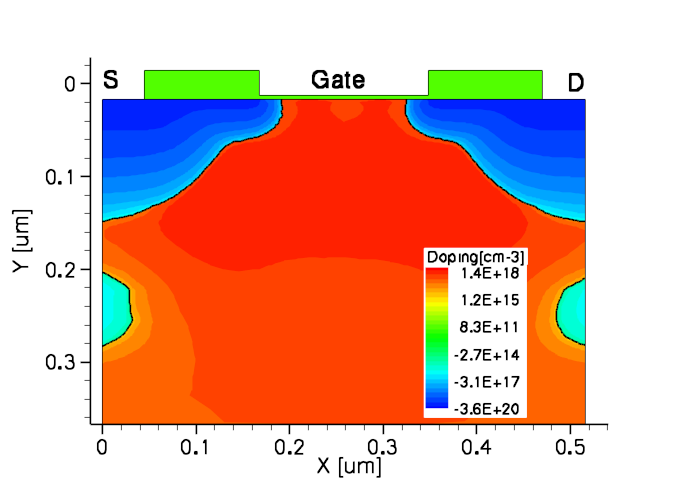}}
\hspace{0.15\textwidth}
\subfloat[][\emph{Mesh}\label{fig: MOS mesh}]
{\includegraphics[width = 0.3\textwidth , height=4.5cm]{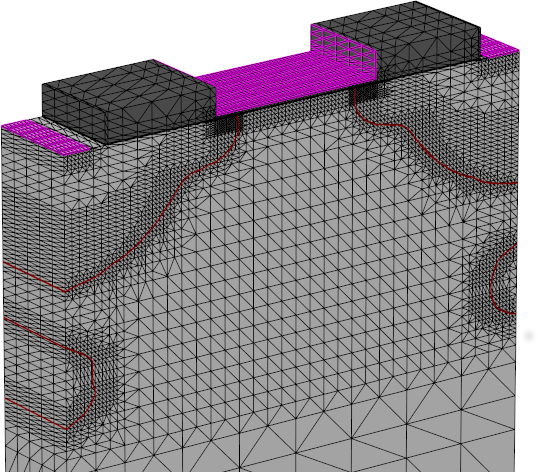}}
\caption{$p$-MOSFET structure: Left: Doping. Right: Mesh.}
\label{fig: pMOSFET}
\end{figure}

\begin{figure}[!h]
\centering
\subfloat[][\emph{$I_D-V_G$ characteristic with $V_D=-0.1V$}\label{fig: IdVg}]
{\includegraphics[width = 0.4\textwidth , height=4.5cm]{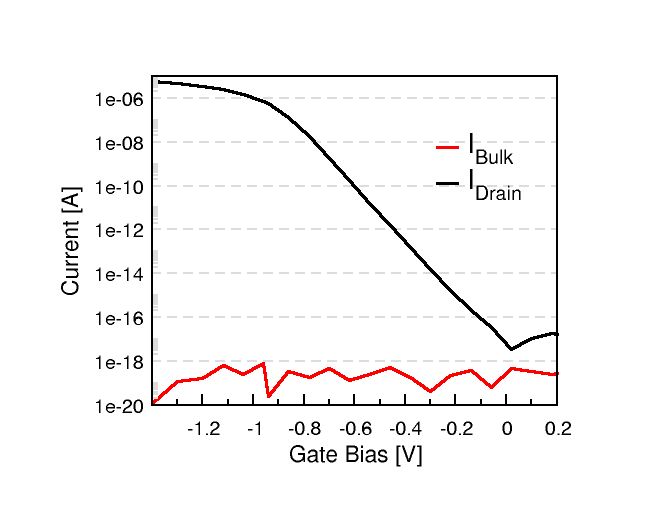}}
\hspace{0.15\textwidth}
\subfloat[][\emph{Method DDFE}\label{fig: Jp DD mos}]
{\includegraphics[width = 0.3\textwidth , height=4.5cm]{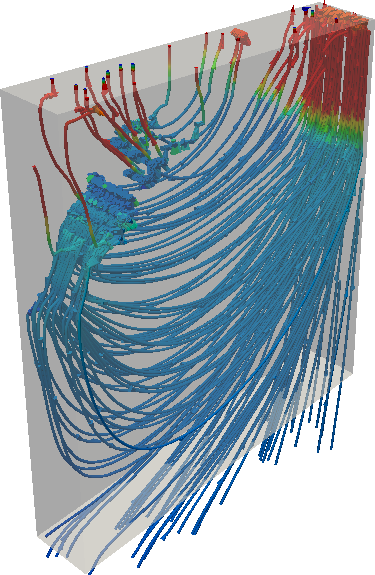}}
\hspace{0.02\textwidth}
{\includegraphics[width = 0.1\textwidth , height=4.5cm]{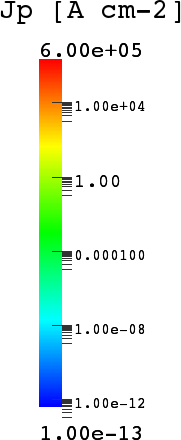}}
\hspace{0.05\textwidth}
\subfloat[][\emph{Method A}\label{fig: Jp SG mos}]
{\includegraphics[width = 0.3\textwidth , height=4.5cm]{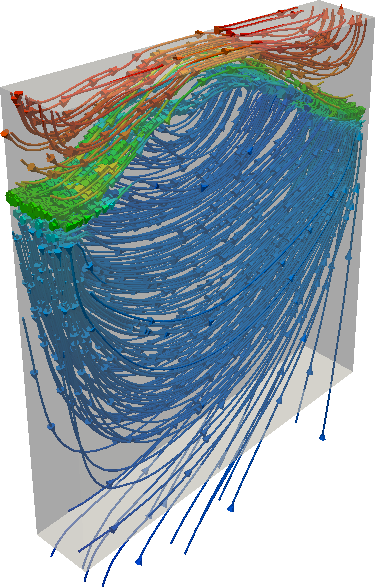}}
\hspace{0.02\textwidth}
{\includegraphics[width = 0.1\textwidth , height=4.5cm]{/img/MOS/JP_LEGENDA}}
\hspace{0.05\textwidth}
\subfloat[][\emph{Method B}\label{fig: Jp UP mos}]
{\includegraphics[width = 0.3\textwidth , height=4.5cm]{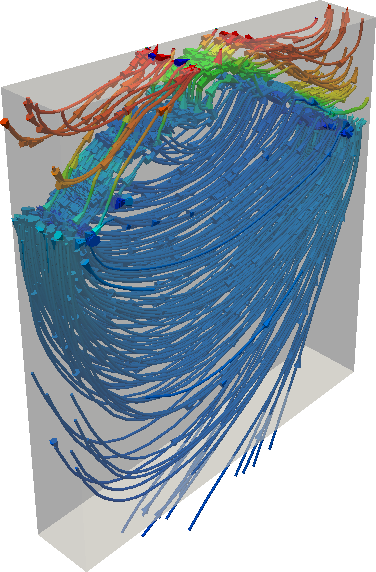}}
\hspace{0.02\textwidth}
{\includegraphics[width = 0.1\textwidth , height=4.5cm]{/img/MOS/JP_LEGENDA}}
\caption{$I_D-V_G$ and hole current density calculation of $p$-MOSFET 
at $V_G=-1.0V$ and $V_D=-0.1V$: Top left: $I_D-V_G$. 
Top right: Method DDFE. Bottom left: Method A. Bottom right: Method B.}
\label{fig: pMOSFET current methods}
\end{figure}

The comparison of the different current computational methods has also been 
carried out on a $p$-channel MOSFET. The doping profiles have been obtained by using a 2-D process simulator with implantation and diffusion steps~\cite{SDeviceManual} with the purpose to have a realistic doping as reported in~\reffig{fig: MOS doping} with  a $Gate Length: 180 nm$ and a 
$Gate Ox: 4.5 nm$. The presence of floating non-compensated $p$-type regions 
in the channel body increases the computational difficulties. The 2D doping profile has been then extruded in three spatial dimensions, and the generated mesh is shown in \reffig{fig: MOS mesh} where the device contact has been highlighted with purple color (the body contact is not shown in the picture but is located in the down face of the silicon region). The Gate contact is then negatively biased to build the standard $I_D-V_G$ with the drain bias kept at 0.1 V. The calculation of the current at the contacts is carried out
by extending to the 3D case the approach proposed 
in~\cite{ContactCurrentRM} for the 2D case. 
The novel formulation, known as \textit{residual method}, 
is based on the approach proposed in~\cite{GalerkMethConsHughes}. 
Computed currents are reported in \reffig{fig: IdVg} (lines) and compared with the commercial tool results~\cite{SDeviceManual} (symbols). The two computed curves are virtually indistinguishable. 
Regarding the calculation of $\vect{J}_p$, 
the numerical difficulties found with Method DDFE are still confirmed as clearly depicted in \reffig{fig: Jp DD mos} not only inside the $p$-type region but also around the floating regions present in the body (the visualization are referred to the bias conditions $V_G=-1.0V$ and $V_D=-0.1V$). The marginality found using formula (\ref{eq:  Jp 3D upwin SG}) in Sect.~\ref{diode} is increased as reported in \reffig{fig: Jp UP mos}: this comes by the fact that the evaluation of the coefficient in the formula (\ref{eq: SG stabilized}) is again undergoing numerical problems related to roundoff error. However Method B is giving a much better current density evaluation with respect to the pure application of the Drift-Diffusion approach at element wise level.
We conclude this section by noting that, again, the best description of the expected physical behavior of the device is obtained by adopting Method A, that turns out to provide an accurate and stable 3D extension of the 1D SG formula, as clearly demonstrated by \reffig{fig: Jp SG mos}.

\subsection{$n$-MOS structure with II}
\begin{figure}[!h]
\subfloat[][\emph{Mesh and Doping}\label{fig: mos geometry 1}]
{\includegraphics[width=0.34\textwidth ,height=4.6cm]
{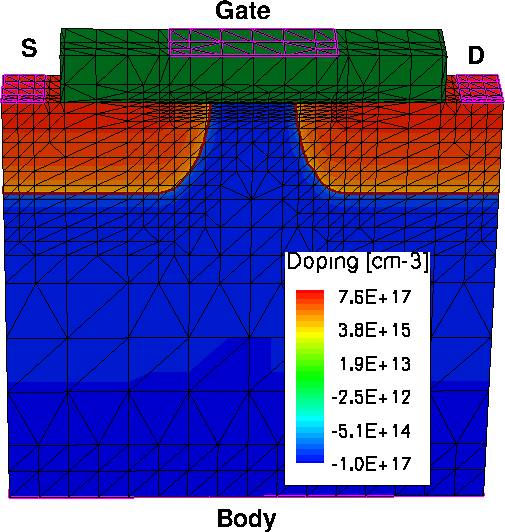}}
\hspace{0.18\textwidth}
\subfloat[][\emph{$I_D-V_G$ characteristics}\label{fig: mos geometry 2}]
{\includegraphics[width = 0.45\textwidth , height=5cm]{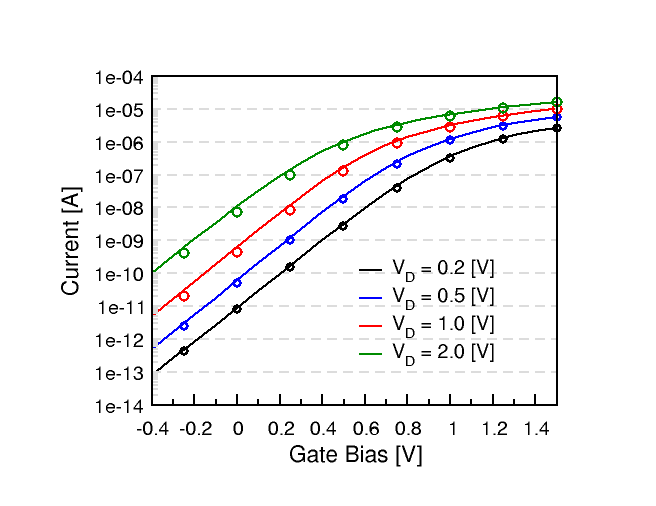}}
\caption{$n$-MOS structure and numerical validation in on-state: Left: Doping and Mesh. Right: $I_D-V_G$.}
\label{fig: mos geometry}
\end{figure}

\begin{figure}[!h]
\centering
\subfloat[][\emph{linear scale}\label{fig: MOS iv_lin}]
{\includegraphics[width = 0.49\textwidth , height=5.5cm]
{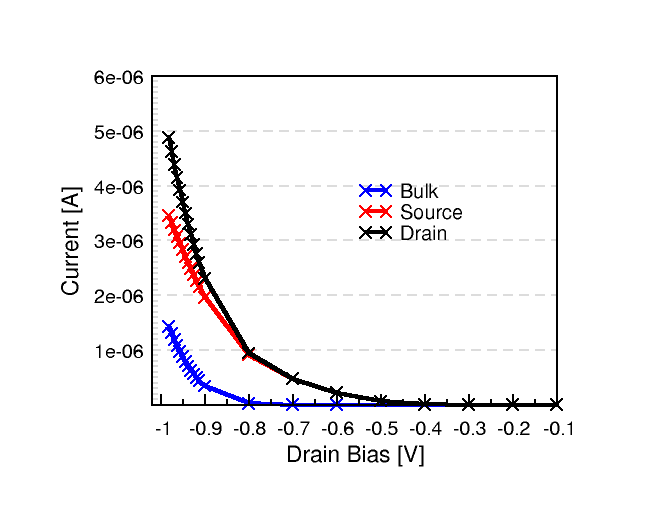}}
\subfloat[][\emph{log-scale}\label{fig: MOS iv_log}]
{\includegraphics[width = 0.49\textwidth , height=5.5cm]
{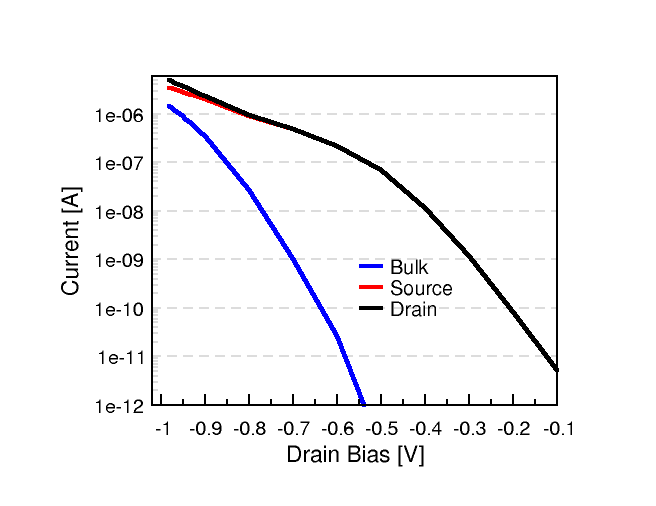}}
\caption{n-MOSFET: off-state characteristic. Left: linear scale.
Right: log-scale.}
\label{fig: nMOSFET_IV_curves}
\end{figure}

\begin{figure}[!h]
\centering
\subfloat[][\emph{$II$, $V_D=-0.5 \, \unit{V}$}\label{fig: MOS II_1}]
{\includegraphics[width = 0.49\textwidth , height=5.5cm]
{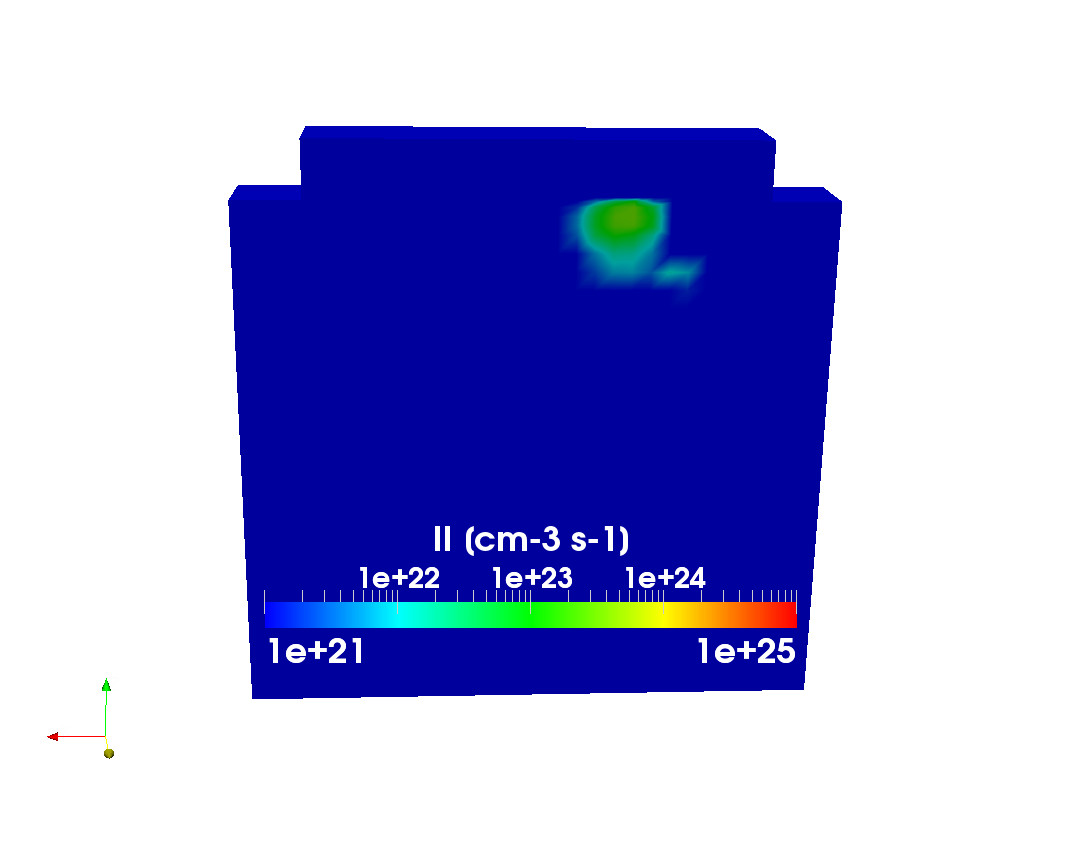}}
\subfloat[][\emph{$II$, $V_D=-0.85 \, \unit{V}$}\label{fig: MOS II_2}]
{\includegraphics[width = 0.49\textwidth , height=5.5cm]
{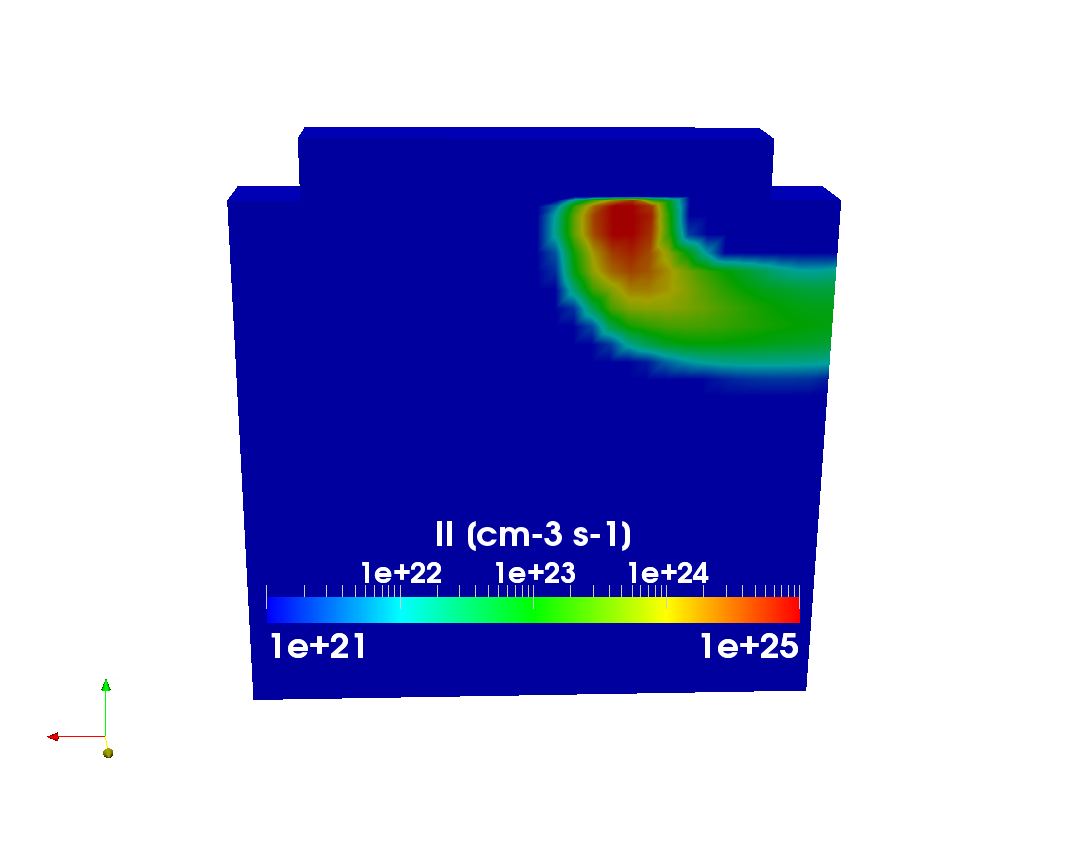}} \\
\subfloat[][\emph{$II$, $V_D=-0.93 \, \unit{V}$}\label{fig: MOS II_3}]
{\includegraphics[width = 0.49\textwidth , height=5.5cm]
{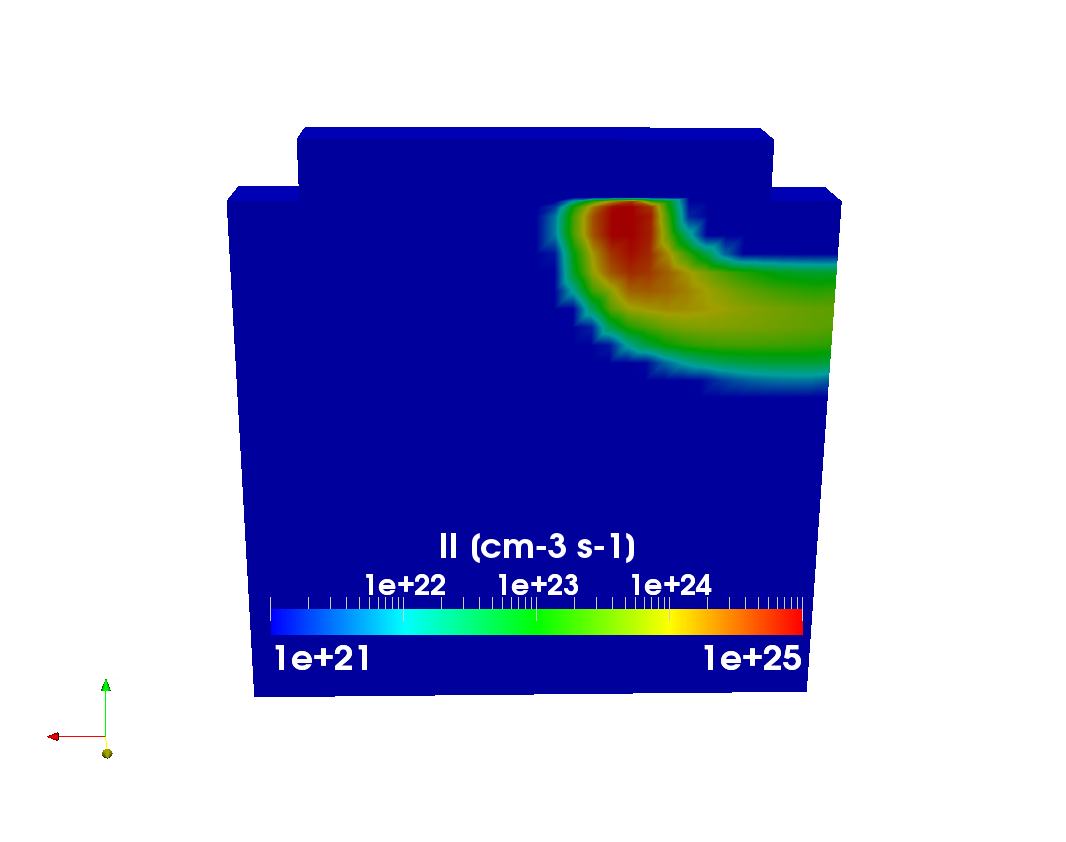}}
\subfloat[][\emph{$II$, $V_D=-0.98 \, \unit{V}$}\label{fig: MOS II_4}]
{\includegraphics[width = 0.49\textwidth , height=5.5cm]
{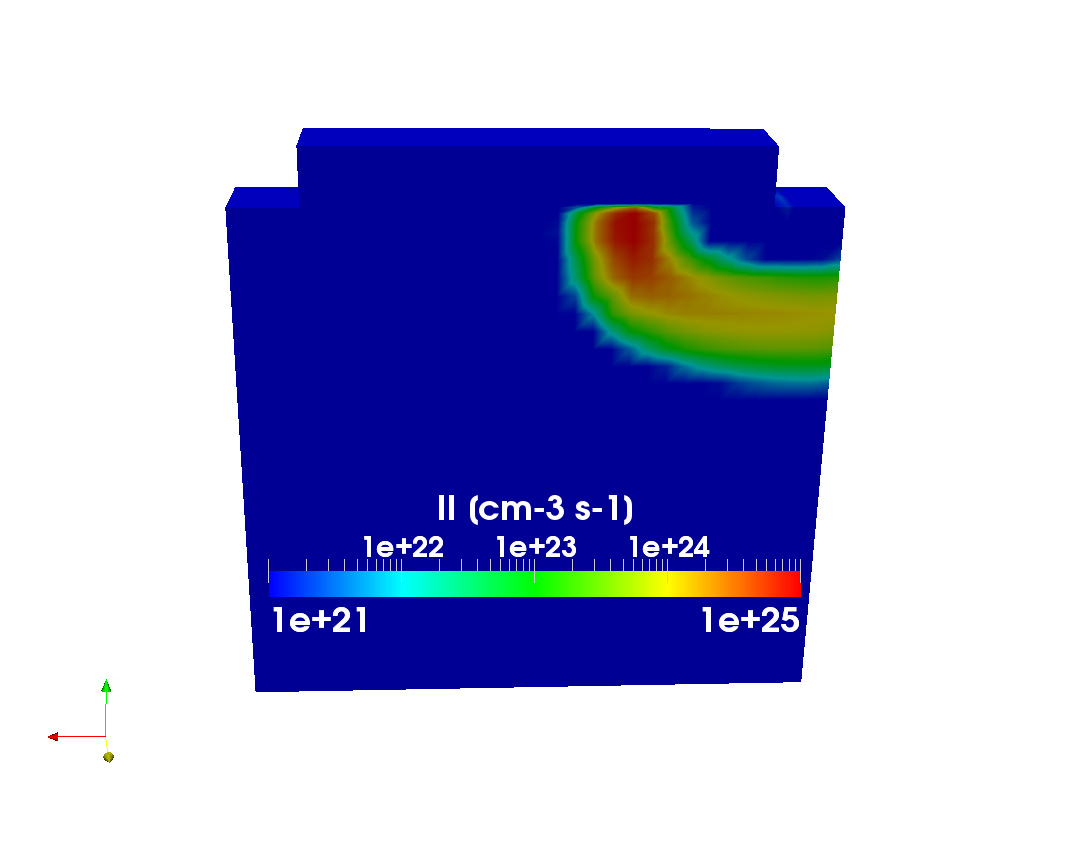}}
\caption{n-MOSFET: II generation term. Top left: $V_D=-0.5 \, \unit{V}$.
Top right: $V_D=-0.85 \, \unit{V}$. Bottom left: $V_D=-0.93 \, \unit{V}$.
Bottom right: $V_D=-0.98 \, \unit{V}$.}
\label{fig: n_mos II}
\end{figure}

\begin{figure}[!h]
\centering
\subfloat[][\emph{$\vect{J}_n$ on-state}\label{fig: MOS J_n onstate}]
{\includegraphics[width = 0.49\textwidth , height=5.5cm]
{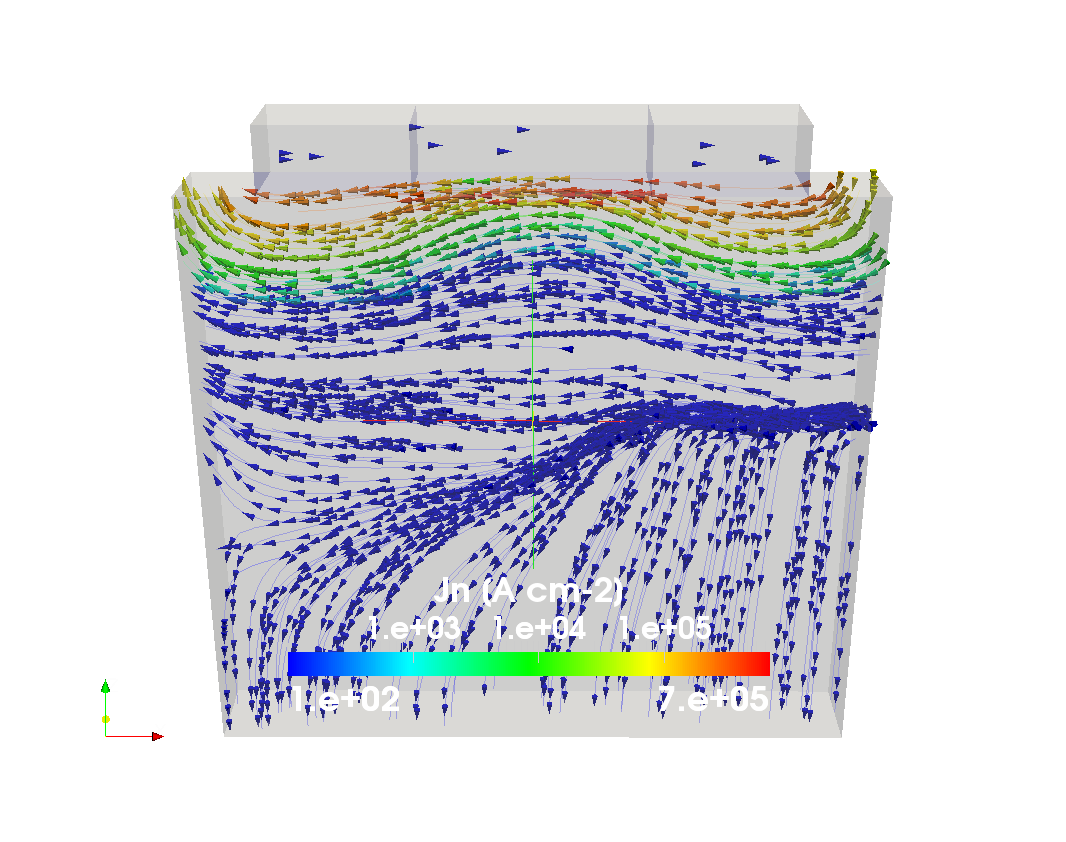}}
\subfloat[][\emph{$\vect{J}_n$ off-state}\label{fig: MOS J_n offstate}]
{\includegraphics[width = 0.49\textwidth , height=5.5cm]
{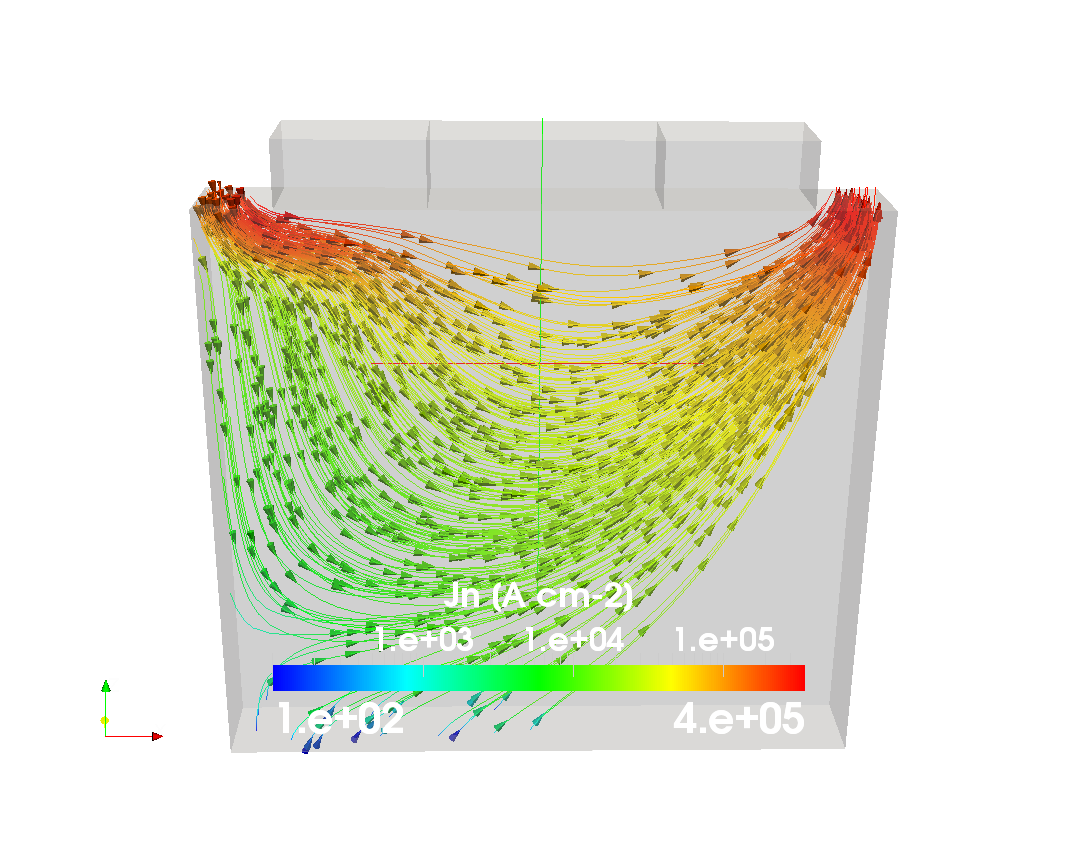}}
\caption{n-MOSFET: electron current density streamlines obtained 
with method A. Left: on-state. Right: off-state.}
\label{fig: nMOSFET_streamlines}
\end{figure}

As a final simulation test we have adopted the $n$-channel MOSFET structure of \reffig{fig: mos geometry 1} ($Gate Length:100 nm; Gate Ox:30 nm$ and contacts $S=source,Gate,D=drain,Body$) with analytical doping profiles. 
Numerical validation of the on-state is reported in \reffig{fig: mos geometry 2} for different values of drain bias. Our results (continuous lines) are compared with those of a commercial tool \cite{SDeviceManual} (symbols) obtaining a very good agreement: electron and hole mobilities are assumed to be constant $\mu_n=1417,\mu_p=470.5 [cm^{2}V^{-1}s^{-1}]$.
\reffig{fig: nMOSFET_IV_curves} shows the results of the simulation in the case  of reverse bias, where the II mechanism is activated as a generation term by employing Method A for the current calculation inside the Gummel map algorithm. The drain-to-bulk breakdown is found to occurr at around $V_D=-1V$ and in the linear scale plot we have highlighted with the cross-symbol the reduction of the drain bias incremental step automatically determined due to high difficulty in convergence when the II generation rate significantly increases and dominates over the other R/G phenomena. \reffig{fig: n_mos II} visualizes the increase 
of the II generation term with the drain inverse voltage: breakdown is starting at the drain-channel junction under the gate and proceeds towards the bulk along the junction profile. 
\reffig{fig: nMOSFET_streamlines} reports the visualization of 
$\vect{J}_n$ in both the off and on states of the device, calculated using Method A: as expected in the on-state the electron current density is confined under the gate from drain-to-source while in the off-state it is flowing from source-to-drain and from source-to-body in accordance to physical insight.

\section{Conclusions}
\label{sec:conclusion}

In this article we have addressed the problem of representing in 3D 
the carrier current density by 
extending the beneficial properties emanating from the classic 1D 
Scharfetter-Gummel difference formula.
To this purpose, we have adopted the Galerkin Finite Element Method 
for the numerical simulation of the Drift-Diffusion model in the 3D
case, and we have proposed two novel methods for
current density evaluation,  denoted method A and B, to which, for 
comparison, we have added also the basic method DDFE using the DD formula.
The three schemes compute a piecewise constant approximation
of the current density over a tetrahedral partition of the device domain.

Method DDFE turns out
to provide the worst results in the test experiments. Method B is a
3D extension of the method of optimal artificial diffusion and gives 
reasonably accurate results. Method A is a natural extension of the 1D SG
approach based on a primal-mixed formulation endowed with a suitable
quadrature formula for the approximate evaluation of the averaged electrical
conductivity. It is by far the best of the three considered approaches, providing
simulation results in excellent agreement with a commercial software.
The three FE formulations are numerically verified in the study of 
realistic 3D structures, also including the presence of Impact Ionization phenomena. Even in this latter case, Method A is able to correctly describe the complex carrier flow patterns inside the device bulk and to track the I-V curve
of the device up to the avalanche breakdown region.

Despite the proposed formulations are illustrated and validated in
the study of the classic DD transport model in semiconductors,
they can be applied in a straightforward manner to the numerical 
treatment of general conservation laws for advective-diffusive fluxes 
where the advective term is in gradient form, as is the case for 
ion electrodiffusion in electrochemistry and biology 
with the Poisson-Nernst-Planck model~\cite{Rubinstein1990} and
hydrodynamic and quantum-corrected charge transport in
semiconductors~\cite{Forghieri1988,deFalcoSacco2005,deFalcoSacco2009}.

\bibliographystyle{plain} 
\bibliography{Bibliografia}

\begin{thebibliography}{10}

\bibitem{SDeviceManual}
{\em Sentaurus Device User Guide}.
\newblock Synopsis Inc., 2013.

\bibitem{Allen:SGformula}
D.~N.~G. Allen and R.~V. Southwell.
\newblock Relaxation methods applied to determine the motion, in two
  dimensions, of a viscous fluid past a fixed cylinder.
\newblock {\em Quart. J. Mech. Appl. Math.}, 8:129--145, 1955.

\bibitem{Bank:Upwinding}
R.~E. Bank, J.~F. Burgler, W.~Fichtner, and R.~K. Smith.
\newblock Some upwinding techniques for finite element approximations of
  convection-diffusion equations.
\newblock {\em Numer. Math.}, 58:185--202, 1990.

\bibitem{BankRoseFichtnersiam1983}
R.E. Bank, D.~Rose, and W.~Fichtner.
\newblock Numerical methods for semiconductor device simulation.
\newblock {\em SIAM Journal on Scientific and Statistical Computing},
  4(3):416--435, 1983.

\bibitem{BankRoseFichtnerieee1983}
R.E. Bank, D.~Rose, and W.~Fichtner.
\newblock Numerical methods for semiconductor device simulation.
\newblock {\em Electron Devices, IEEE Transactions on}, 30(9):1031--1041, 1983.

\bibitem{Brezzicmame1989}
F.~Brezzi, L.~Marini, and P.~Pietra.
\newblock Numerical simulation of semiconductor devices.
\newblock {\em Comp. Meths. Appl. Mech. Engrg}, 75:493--514, 1989.

\bibitem{Brezzi1989}
F.~Brezzi, L.~Marini, and P.~Pietra.
\newblock Two-dimensional exponential fitting and applications to
  drift-diffusion models.
\newblock {\em SIAM Journal on Numerical Analysis}, 26(6):1342--1355, 1989.

\bibitem{BrezziLectNotes1991}
F.~Brezzi, L.D. Marini, P.A. Markowich, and P.~Pietra.
\newblock On some numerical problems in semiconductor device simulation.
\newblock In G.~Toscani, V.~Boffi, and S.~Rionero, editors, {\em Mathematical
  Aspects of Fluid and Plasma Dynamics}, volume 1460 of {\em Lecture Notes in
  Mathematics}, pages 31--42, 1991.

\bibitem{Brooks1982}
A.~N. Brooks and T.~J.R. Hughes.
\newblock Streamline {U}pwind/{P}etrov-{G}alerkin formulations for convection
  dominated flows with particular emphasis on the incompressible
  {N}avier-{S}tokes equations.
\newblock {\em Computer Methods in Applied Mechanics and Engineering},
  32(1\u20133):199 -- 259, 1982.

\bibitem{ITRS2013}
P.~Cogez, M.~Graef, B.~Huizing, R.~Mahnkopf, H.~Ishiuchi, N.~Ikumi, S.~Choi,
  J.~H. Choi, C.~H. Diaz, Y.~C. See, P.~Gargini, T.~Kingscott, and I.~Steff.
\newblock The {I}nternational {T}echnology {R}oadmap for {S}emiconductors:
  2013.
\newblock Technical report, Jointly sponsored by ESIA, JEITA, KSIA, TSIA and
  SIA, 2013.

\bibitem{deFalcoSacco2009}
C.~de~Falco, J.~W. Jerome, and R.~Sacco.
\newblock Quantum-corrected drift-diffusion models: Solution fixed point map
  and finite element approximation.
\newblock {\em Journal of Computational Physics}, 228(5):1770 -- 1789, 2009.

\bibitem{deFalcoSacco2005}
C.~de~Falco, A.~L. Lacaita, E.~Gatti, and R.~Sacco.
\newblock Quantum-corrected drift-diffusion models for transport in
  semiconductor devices.
\newblock {\em J. Comp. Phys.}, 204:533--561, 2005.

\bibitem{Forghieri1988}
A.~Forghieri, R.~Guerrieri, P.~Ciampolini, A.~Gnudi, M.~Rudan, and
  G.~Baccarani.
\newblock A new discretization strategy of the semiconductor equations
  comprising momentum and energy balance.
\newblock {\em Computer-Aided Design of Integrated Circuits and Systems, IEEE
  Transactions on}, 7(2):231--242, Feb 1988.

\bibitem{GummelMap}
H.~K. Gummel.
\newblock A self-consistent iterative scheme for one-dimensional steady state
  transistor calculations.
\newblock {\em IEEE Trans. Electron Devices}, pages 455--465, 1964.

\bibitem{Gummel:SignAnalys}
H.~K. Gummel and D.~Scharfetter.
\newblock Large-signal analysis of a silicon read diode oscillator.
\newblock {\em IEEE Trans. Electron Devices}, pages 64--77, 1969.

\bibitem{ContactCurrentRM}
R.~Gusmeroli and A.~S. Spinelli.
\newblock Accurate boundary integrals calculation in semiconductor device
  simulation.
\newblock {\em IEEE Trans. Electron Devices}, 53:1730--1733, 2006.

\bibitem{GalerkMethConsHughes}
T.J.R. Hughes, G.~Engel, L.~Mazzei, and M.G. Larson.
\newblock The continuous {G}alerkin method is locally conservative.
\newblock {\em Journal of Computational Physics}, 163:467--488, 2000.

\bibitem{Jerome:AnalyCharTran}
J.~W. Jerome.
\newblock {\em Analysis of Charge Transport}.
\newblock Springer, 1996.

\bibitem{markowich1986stationary}
P.A. Markowich.
\newblock {\em The Stationary Semiconductor Device Equations}.
\newblock Computational Microelectronics. Springer-Verlag, 1986.

\bibitem{Mauri2014}
A.~Mauri, R.~Sacco, and M.~Verri.
\newblock Electro-thermo-chemical computational models for {3D} heterogeneous
  semiconductor device simulation.
\newblock arXiv:1307.3096, to appear in Applied Mathematical Modeling.

\bibitem{Micheletti1995}
S.~Micheletti, A.~Quarteroni, and R.~Sacco.
\newblock Current-voltage characteristics simulation of semiconductor devices
  using domain decomposition.
\newblock {\em Journal of Computational Physics}, 119(1):46 -- 61, 1995.

\bibitem{NovielliIEDM2013}
G.~Novielli, A.~Ghetti, E.~Varesi, A.~Mauri, and R.~Sacco.
\newblock Atomic migration in phase change materials.
\newblock In {\em Electron Devices Meeting (IEDM), 2013 IEEE International},
  pages 22.3.1--22.3.4, Dec 2013.

\bibitem{quarteroni:NumApprox}
A.~Quarteroni and A.~Valli.
\newblock {\em Numerical Approximation of Partial Differential Equations}.
\newblock Springer, 2008.

\bibitem{Roberts:MixedHybrid}
J.~Roberts and J.~M. Thomas.
\newblock Mixed and hybrid methods.
\newblock {\em Handbook of Numerical Analysis}, 2:523--633, 1991.

\bibitem{Rubinstein1990}
I.~Rubinstein.
\newblock {\em Electrodiffusion of Ions}.
\newblock SIAM, 1990.

\bibitem{Sacco2014ABB}
R.~Sacco, L.~Carichino, C.~de~Falco, M.~Verri, F.~Agostini, and T.~Gradinger.
\newblock A multiscale thermo-fluid computational model for a two-phase cooling
  system.
\newblock {\em Computer Methods in Applied Mechanics and Engineering},
  282(0):239 -- 268, 2014.

\bibitem{selberherr:SimSem}
S.~Selberherr.
\newblock {\em Analysis and Simulation of Semiconductor Devices}.
\newblock Springer-Verlag, Wien, 1984.

\bibitem{Zikatanov:EAFE1}
J.~Xu and L.~Zikatanov.
\newblock A monotone finite element scheme for convection-diffusion equations.
\newblock {\em Mathematics of Computation}, 68:1429--1446, 1999.

\bibitem{Zikatanov:EAFE2}
L.~T. Zikatanov and R.~D. Lazarov.
\newblock An exponential fitting scheme for general convection-diffusion
  equations on tetrahedral meshes.
\newblock {\em Computational and Applied Mathematics}, 1:60--69, 2012.

\end{thebibliography}

\end{document}